\documentclass{elsart}

\usepackage{amssymb,amsmath}

\newtheorem{theorem}{Theorem}[section]

\newtheorem{proposition}[theorem]{Proposition}
\newtheorem{corollary}[theorem]{Corollary}

\newtheorem{question}[theorem]{Question}

\journal{Acta Mathematica Hungarica}
\begin{document}

\begin{frontmatter}

\title{Metrizable compacta in the space of continuous functions with the
topology of pointwise convergence}

\author{V.V. Mykhaylyuk}

\address{Chernivtsi National University, Department of Mathematical Analysis,
Kotsjubyns'koho 2, Chernivtsi 58012, Ukraine\\e-mail:
mathan@ukr.net}

\begin{abstract}
We prove that every point-finite family of nonempty functionally
open sets in a topological space $X$ has the cardinality at most
an infinite cardinal $\kappa$ if and only if $w(X)\leq\kappa$ for
every Valdivia compact space $Y\subseteq C_p(X)$. Correspondingly
a Valdivia compact space $Y$ has the weight at most an infinite
cardinal $\kappa$ if and only if every point-finite family of
nonempty open sets in $C_p(Y)$ has the cardinality at most
$\kappa$, that is $p(C_p(Y))\leq \kappa$. Besides, it was proved
that $w(Y)=p(C_p(Y))$ for every linearly ordered compact $Y$. In
particular, a Valdivia compact space or linearly ordered compact
space $Y$ is metrizable if and only if $p(C_p(Y))=\aleph_0$. This
gives answer to a question of O.~Okunev and V.~Tkachuk.

\end{abstract}

\end{frontmatter}

\section{Introduction}
\label{Introduction}

During the study of the Baire classification of separately
continuous functions $f:X\times Y\to \mathbb R$, when $Y$
satisfies compactness type conditions, it was established that
spaces $X$ which satisfy a condition of a countable chain type,
play an important role (see [1,2,3]). The Baire classification of
these mappings was obtained replacing of $Y$ with its continuous
metrizable image in the space $C_p(X)$ of all continuous functions
on $X$ with the pointwise convergence topology.

This shows that it is natural to try to find necessary and
sufficient conditions on a topological space $X$ for metrizability
of any compact subspace of $C_p(X)$ or any compact subspace of
$C_p(X)$ which belongs to a given class such as Eberlein compact
spaces, Corson compact spaces or Valdivia compact spaces.

A compact space $X$ which is homeomorphic to a weakly compact
subset of a Banach space is called {\it an Eberlein compact
space}. A compact space $X$ is called {\it a Corson compact space}
if it is homeomorphic to a compact space $Z\subseteq {\mathbb
R}^T$ such that $|{\rm supp}\,z|\leq \aleph_0$ for any $z\in Z$.
By ${\rm supp}\,z$ we denote the support $\{t\in T:z(t)\ne 0\}$ of
a function $z:T\to\mathbb R$. A compact space $X$ is called {\it a
Valdivia compact space} if it is homeomorphic to a compact space
$Z\subseteq {\mathbb R}^T$ such that the set $\{z\in Z:|{\rm
supp}\,z|\leq \aleph_0\}$ is dense in $Z$.

Analogous questions were investigated in [4,5,6,7],where
connections between cardinal properties of a completely regular
space $X$ and the space $C_p(X)$ were studied. The point-finite
cellularity $p(X)= \sup \{ |\tau|: \tau \mbox{\,\, is\,\,a
\,\,point-finite\,\,family}$
$\mbox{\,\,of\,\,nonempty\,\,open\,\,sets\,\,in\,\,}X\}$ of a
topological space $X$, which was introduced in [4], is important
in these papers. In particular, the following questions were posed
in [7, questions 5.2, 5.4]:

\begin{question}
{\it Let $X$ be a Souslin continuum, i.e. a non-separable linearly
ordered perfectly normal compact space. It is true that
$p(C_p(X))=\aleph_0$?}
\end{question}

\begin{question}
{\it Let $X$ be a Corson compact space  such that
$p(C_p(X))=\aleph_0$. Must $X$ be metrizable?}
\end{question}

Question 1.2 was answered in [8]. It is proved in [8] that the
weight of a Corson compact space $X$ is equal to the point-finite
cellularity of $C_p(X)$.

In this paper we study necessary and sufficient conditions on a
topological space $X$ for metrizability of any Valdivia compact
space $Y\subseteq C_p(X)$. Using a technique of a dependence of
mappings upon a certain number of coordinates we prove that any
point-finite family of functionally open sets in a topological
space $X$ has the cardinality at most an infinite cardinal
$\kappa$ if and only if each Valdivia compact space $Y\subseteq
C_p(X)$ has the topological weight at most $\kappa$. In
particular, $p(X)=\aleph_0$ for a completely regular space $X$ if
and only if each Valdivia compact space $Y\subseteq C_p(X)$ is
metrizable. Passing to a dual problem on study of connections
between properties of a compact space $Y$ and the space $C_p(Y)$,
we obtain that a Valdivia compact space $Y$ has the weight at most
an infinite cardinal $\kappa$ if and only if $p(C_p(Y))\leq
\kappa$. It implies immediately the following fact: a Valdivia
compact space $Y$ is metrizable if and only if
$p(C_p(Y))=\aleph_0$. Besides, it was proved that $w(Y)=p(C_p(Y))$
for every linearly ordered compact space $Y$ (this yields a
negative answer to Question 5.2 from [7]).

\section{A dependence of mappings upon a certain number of coordinates}

In this section we prove a theorem which will be a technical tool
in the proof of the main result and based on a notion of  a
dependence of mappings upon a certain number of coordinates.

Let $\kappa$ be an infinite cardinal number, and let $Z$, $T$ be
sets, $Y\subseteq {\mathbb R}^T$ and $f:Y\to Z$. We say that {\it
$f$ is concentrated on a set $S\subseteq T$} if for any $y',
y''\in Y$ the condition $y'|_S = y''|_S$ implies $f(y')=f(y'')$.
Besides, if $|S|\leq \kappa$ then {\it $f$ depends upon $\kappa$
coordinates.}

Let $X$ be a set and $g:X\times Y\to Z$. Then {\it $g$ is
concentrated on a set $S\subseteq T$ in the second variable} if
$g(x,y')=g(x,y'')$ for any $x\in X$ and $y',y''\in Y$ with $y'|_S
= y''|_S$; and {\it $g$ depends upon $\kappa$ coordinates in the
second variable} if $|S|\leq \kappa$ for some $S$.

A subset $A$ of a topological space $X$ is called {\it
functionally open} if there exists a continuous function $f:X\to
[0,1]$ such that $A=f^{-1}((0,1])$.

The following result was obtained in [9] for $\kappa=\aleph_0$.

\begin{theorem}
{\it Let $X$ be a topological space in which every point-finite
family of nonempty functionally open sets has the cardinality at
most an infinite cardinal number $\kappa$, let $Y\subseteq
{\mathbb R}^T$ be a compact space, $B=\{y\in Y: |{\rm
supp}\,y|\leq \kappa\}$ and $Y=\overline{B}$. Then every
separately continuous function $f:X\times Y\to\mathbb R$ depends
upon $\kappa$ coordinates in the second variable.}
\end{theorem}

{\bf Proof.} First we prove that for every $\varepsilon
>0$ there exists a set $R_{\varepsilon}\subseteq T$ such that
$|R_{\varepsilon}|\leq \kappa$ and for any $b',b''\in B$ the
condition $b'|_{R_{\varepsilon}}=b''|_{R_{\varepsilon}}$ implies
$|f(x,b')-f(x,b'')|\leq \varepsilon$ for every $x\in X$.

Suppose, contrary to our claim, that there exists an $\varepsilon
>0$ such that for every set $S\subseteq T$ with $|S|\leq
\kappa$ there exist an $x\in X$ and $b',b''\in B$ such that
$b'|_S=b''|_S$ and $|f(x,b')-f(x,b'')|>\varepsilon$.

Denote by $T(\kappa^+)$ the first ordinal of cardinality
$\kappa^+$. Using the transfinite induction we construct a system
$\{S_{\alpha}:\alpha<{T(\kappa^+)}\}$ of sets $S_{\alpha}\subseteq
T$ and sets $\{b_{\alpha}:\alpha<{T(\kappa^+)}\}$,
$\{c_{\alpha}:\alpha<{T(\kappa^+)}\}$ and
$\{x_{\alpha}:\alpha<{T(\kappa^+)}\}$ of points $b_{\alpha},
c_{\alpha}\in B$ and $x_{\alpha}\in X$ such that

$(a)$\,\,\,$|S_{\alpha}|\leq \kappa$ for any $\alpha<T(\kappa^+)$;

$(b)$\,\,\,$b_{\alpha}|_{S_{\alpha}}=c_{\alpha}|_{S_{\alpha}}$ for
any $\alpha<T(\kappa^+)$;

$(c)$\,\,\,$S_{\alpha}\subseteq S_{\beta}$ for any
$\alpha<\beta<T(\kappa^+)$;

$(d)$\,\,\,${\rm supp}\,b_{\alpha}\subseteq S_{\alpha+1}$, ${\rm
supp}\,c_{\alpha}\subseteq S_{\alpha+1}$ for any
$\alpha<T(\kappa^+)$;

$(e)$\,\,\,$|f(x_{\alpha}, b_{\alpha})-f(x_{\alpha},
c_{\alpha})|>\varepsilon$ for any $\alpha <T(\kappa^+)$.

Fix a set $S_1\subseteq T$ with $|S_1|\leq \kappa$. By the
assumption, there exist points $x_1\in X$ and $b_1, c_1\in B$ such
that $b_1|_{S_1}=c_1|_{S_1}$ and
$|f(x_1,b_1)-f(x_1,c_1)|>\varepsilon$. Put $S_2=S_1\cup{\rm
supp}\,b_1\cup{\rm supp}\, c_1$. Clearly  $|S_2|\leq \kappa$.
Choose points $x_2\in X$ and $b_2,c_2\in B$ such that
$b_2|_{S_2}=c_2|_{S_2}$ and $|f(x_2,b_2)-f(x_2,c_2)|>\varepsilon$.

Suppose that for some $\beta<T(\kappa^+)$ the sets
$\{S_{\alpha}:\alpha<\beta\}$, $\{b_{\alpha}:\alpha<\beta\}$,
$\{c_{\alpha}:\alpha<\beta\}$ and $\{x_{\alpha}:\alpha<\beta\}$
are constructed. Put
$S_{\beta}=\bigcup\limits_{\alpha<\beta}(S_{\alpha}\cup{\rm
supp}\,b_{\alpha}\cup{\rm supp}\,c_{\alpha})$. Since for
$\alpha<\beta$ all the sets $S_{\alpha}$, ${\rm supp}\,b_{\alpha}$
and ${\rm supp}\,c_{\alpha}$ have cardinality at most $\kappa$, we
have the inequality $|S_{\beta}|\leq \kappa$. By the assumption,
there exist points $x_{\beta}$ and $b_{\beta}, c_{\beta}\in B$
such that $b_{\beta}|_{S_{\beta}}=c_{\beta}|_{S_{\beta}}$ and
$|f(x_{\beta},b_{\beta})-f(x_{\beta},c_{\beta})|>\varepsilon$.

The continuity of $f$ in variable $x$ and condition $(e)$ imply
that for every $\alpha<T(\kappa^+)$ there exists an functionally
open neighborhood $U_{\alpha}$ of $x_{\alpha}$ in $X$ such that
$|f(x,b_{\alpha})-f(x,c_{\alpha})|>\varepsilon$ for any $x\in
U_{\alpha}$. According to the theorem condition, the family
$(U_{\alpha}:\alpha<\omega)$ is not point-finite. Thus there exist
a point $x_0\in X$ and a strictly increasing sequence
$(\alpha_n)^{\infty}_{n=1}$ of ordinals $\alpha_n<T(\kappa^+)$
such that $|f(x_0,b_{\alpha_n})-f(x_0,c_{\alpha_n})|>\varepsilon$
for any $n\in \mathbb N$.

Put $T_n=S_{\alpha_n}$, $v_n=b_{\alpha_n}$ and $w_n=c_{\alpha_n}$
by $n\in \mathbb N$. Using the compactness of $Y$ and the
continuity of $f^{x_0}:Y\to \mathbb R$, $f^{x_0}(y)=f(x_0,y)$,
choose a finite set $T_0\subseteq T$ such that
$|f(x_0,y')-f(x_0,y'')|<\varepsilon$ for any $y',y''\in Y$ with
$y'|_{T_0}=y''|_{T_0}$. It follows from
$|f(x_0,v_n)-f(x_0,w_n)|>\varepsilon$ that $v_n|_{T_0}\ne
w_n|_{T_0}$. But $v_n|_{T_n}= w_n|_{T_n}$ according to $(b)$ and
$v_n|_{T\setminus T_{n+1}}=w_n|_{T\setminus T_{n+1}}$ according to
$(c)$ and $(d)$. Thus \mbox{ $T_0\cap(T_{n+1}\setminus T_n)\ne\O$}
for any $n\in \mathbb N$. Since the sequence
$(T_n)^{\infty}_{n=1}$ is increasing, $T_0$ is infinite, which is
impossible. Hence the existence of $R_{\varepsilon}$ is proved.

Put $R_0=\bigcup\limits_{n=1}^{\infty}R_{\frac{1}{n}}$. Clearly
$f(x,b')=f(x,b'')$ for any $x\in X$, $b',b''\in B$ with
$b'|_{R_0}=b''|_{R_0}$. Fix points $x\in X$, $y',y''\in Y$ such
that $y'|_{R_0}=y''|_{R_0}$. Since a function $f^x:Y\to\mathbb R$,
$f^x(y)=f(x,y)$ is continuous on the compact space $Y\subseteq
{\mathbb R}^T$, there exists an at most countable set
$T_0\subseteq T$ such that $f(x,y_1)=f(x,y_2)$ for any $y_1,y_2\in
Y$ with $y_1|_{T_0}=y_2|_{T_0}$. Consider a continuous mapping
$\varphi:Y\to {\mathbb R}^{T_0\cup R_0}$, $\varphi(y)=y|_{T_0\cup
R_0}$. Since $Y=\overline{B}$, we have the inequality
$\varphi(Y)=\overline{\varphi(B)}$. Note that $B$ is an
$\kappa$-compact set, that is, the closure $\overline{A}$ in $B$
of any set $A\subseteq B$ with $|A|\leq \kappa$ is a compact set.
Therefore $\varphi(B)$ is an $\kappa$-compact set. It follows from
$|T_0\cup R_0|\leq \kappa$ that $w(\varphi(Y))\leq \kappa$. The
density of $\varphi(B)$ in $\varphi(Y)$ implies the existence of a
dense in $\varphi(Y)$ set $C\subseteq \varphi(B)$ such that
$|C|\leq \kappa$. Then $\varphi(Y)\subseteq \overline{C}\subseteq
\varphi(B)$. Thus $\varphi(Y)=\varphi(B)$. Therefore there exist
points $b',b''\in B$ such that $b'|_{T_0\cup R_0}=y'|_{T_0\cup
R_0}$ and $b''|_{T_0\cup R_0}=y''|_{T_0\cup R_0}$. Since
$b'|_{T_0}=y'|_{T_0}$, $b''|_{T_0}=y''|_{T_0}$ and
$b'|_{R_0}=b''|_{R_0}$, we have $f(x,y')=f(x,b')$,
$f(x,b'')=f(x,y'')$ and $f(x,b')=f(x,b'')$. Hence
$f(x,y')=f(x,y'')$. Thus $f$ is concentrated on $R_0$ and
therefore $f$ depends upon $\kappa$ coordinates in the second
variable.\hfill$\diamondsuit$

\section{Main results}

\begin{theorem}
{\it Let $X$ be a topological space and let $\kappa$ be an
infinite cardinal number. Then the following conditions are
equivalent:

$(i)$ every point-finite family of nonempty functionally open in
$X$ sets has the cardinality at most $\kappa$;

$(ii)$ $w(Y)\leq\kappa$ for any Valdivia compact space $Y\subseteq
C_p(X)$;

$(iii)$ $w(Y)\leq\kappa$ for any Eberlein compact space
$Y\subseteq C_p(X)$.}
\end{theorem}

{\bf Proof.} $(i)\Longrightarrow (ii)$. By the definition, a
Valdivia compact space $Y\subseteq C_p(X)$ is homeomorphic to a
compact space $Z\subseteq {\mathbb R}^T$ such that the set $\{z\in
Z: |{\rm supp}\,z|\leq\aleph_0\}$ is dense in $Z$. Fix a
homeomorphism $\varphi:Z\to Y$. For every $x\in X$ and $z\in Z$ we
put $f(x,z)=\varphi(z)(x)$. Clearly, $f$ is a separately
continuous function on $X\times Z$. By Theorem 2.1, $f$ depends
upon $\kappa$ coordinates in the second variable. This implies
that the map $\psi:Z\to {\mathbb R}^S$, $\psi(z)=z|_S$, is an
injection and hence homeomorphism onto $\tilde{Z}=\psi(Z)$.
Therefore $w(Z)=w(\tilde{Z})\leq\kappa$.

Since every Eberlein compact space is a Valdivia compact space,
the implication $(ii)\Longrightarrow (iii)$ is obvious.

The implication $(iii)\Longrightarrow (i)$ was proved in [5,
Theorem 1].

\begin{corollary}
{\it Let $X$ be a completely regular space. Then the following
conditions are equivalent:

$(i)$ $p(X)=\aleph_0$;

$(ii)$ each Valdivia compact space $Y\subseteq C_p(X)$ is
metrizable.}
\end{corollary}

 A topological space $X$ is called {\it
pseudocompact} if every continuous function $f:X\to \mathbb R$ is
bounded; {\it countably compact} if a finite subcover can be
chosen from any countable open cover of $X$. A topological space
$X$ is called {\it $\sigma$-pseudocompact} if it is the union of a
sequence of pseudocompact subspaces of $X$.

\begin{corollary}
{\it Let $X$ be a completely regular space with $p(X)=\aleph_0$
which consists a $\sigma$-pseudocompact dense subspace. Then every
countably compact space $Y\subseteq C_p(X)$ is metrizable.}
\end{corollary}

{\bf Proof.} By [10, Theorem III.4.23] every countably compact
space $Y\subseteq C_p(X)$ is an Eberlein compact space. It remains
to use Corollary 3.2.\hfill$\diamondsuit$

By ${\mathcal K}(X)$ we denote the system of all compact subsets
of a topological space $X$. A topological space $X$ is called {\it
$\mathcal K$-countably-determined} if there exist a subset $Y$ in
the topological space ${\mathbb N}^{\mathbb N}$ and a set-valued
mapping $F:Y\to {\mathcal K}(X)$ such that for every open in $X$
set $U$ the set $\{y\in Y: F(y)\subseteq U\}$ is open in $Y$ and
$X=\bigcup\limits_{y\in Y}F(y)$.

\begin{corollary}
{\it Let $X$ be a completely regular space with $p(X)=\aleph_0$
which contains a ${\mathcal K}$-countably-determined dense
subspace. Then every compact space $Y\subseteq C_p(X)$ is
metrizable.}
\end{corollary}

{\bf Proof.} Fix a dense in $X$ ${\mathcal
K}$-countably-determined subspace $\tilde{X}$ and a compact space
$Y\subseteq C_p(X)$. Consider the continuous identity mapping
$\varphi: C_p(X)\to C_p(\tilde{X})$. Put $\tilde{Y}=\varphi(Y)$.
It follows from density $\tilde{X}$ in $X$ that $\varphi|_Y$ is a
bijection and by [11, Theorem 3.1.13] $\varphi$ is a
homeomorphism. By [12, Theorem 3.7] $\tilde{Y}$ is a Corson
compact space. Thus $Y$ is a Corson compact space and by Corollary
3.2 it is metrizable.\hfill$\diamondsuit$

Let $D$ be a discrete space of cardinality $\kappa>\aleph_0$. By
$L_{\kappa}$ we will denote an one-point Lindel\"{o}ffication of
$D$, that is $L_{\kappa}=D\bigcup \{\infty\}$ and all
neighborhoods of $\infty$ in $L_{\kappa}$ are the complements in
$L_{\kappa}$ to all at most countable subsets of $D$, besides the
discrete topology on $D$ coincides with the topology on $D$
induced by $L_{\kappa}$. A topological space $X$ is called {\it
primary Lindel\"{o}ff} if $X$ is a continuous image of a closed
subset of the topological product ${L_{\kappa}}^{\mathbb N}$ for
some $\kappa>\aleph_0$. Note that by [10, Corollary IV.3.17] for a
primary Lindel\"{o}ff space $X$ every compact space $Y\subseteq
C_p(X)$ is a Corson compact space.

The following result can be obtained analogously.

\begin{corollary}
 {\it Let $X$ be a completely regular space with
$p(X)=\aleph_0$ which contains a primary Lindel\"{o}ff dense
subspace. Then every compact space $Y\subseteq C_p(X)$ is
metrizable.}
\end{corollary}

Recall that a set $F$ of functions defined on a set $X$ {\it
separates points on $X$} if for arbitrary distinct points
$x',x''\in X$ there exists a function $f\in F$ such that $f(x')\ne
f(x'')$.

\begin{theorem}
{\it Let $Y$ be a Valdivia compact space and let $\kappa$ be an
infinite cardinal number. Then the following conditions are
equivalent:

$(i)$ $w(Y)\leq \kappa$;

$(ii)$ $p(C_p(Y))\leq\kappa$;

$(iii)$ there exists a space $X\subseteq C_p(Y)$ with
$p(X)\leq\kappa$ which separates points on $Y$.}
\end{theorem}

{\bf Proof.} $(i)\Longrightarrow (ii)$. Suppose that
$w(Y)\leq\kappa$. If $Y$ is infinite then $d(C_p(Y))\leq
w(Y)\leq\kappa$. If $Y$ is finite then
$d(C_p(Y))=\aleph_0\leq\kappa$ anyway. Thus $d(C_p(Y))\leq\kappa$.
Therefore $p(C_p(Y))\leq\kappa$.

The implication $(ii)\Longrightarrow (iii)$ is obvious.

$(iii)\Longrightarrow (i)$. Fix a subspace $X\subseteq C_p(Y)$
with $p(X)\leq\kappa$ which separates points on $Y$. Consider a
continuous mapping $\varphi:Y\to C_p(X)$, $\varphi(y)(x)=x(y)$.
Since $X$ separates points on $Y$, the map
$\varphi:Y\to\varphi(Y)$ is a bijection and by [11, Theorem
3.1.13] $\varphi$ is a homeomorphism. Thus $\varphi(Y)$ is an
Valdivia compact space and $w(Y)\leq \kappa$ by Theorem
3.1.\hfill$\diamondsuit$

Note that Theorem 3.6 is not true for an arbitrary compact space
$Y$. Indeed, $2^{\aleph_0}=w(Y)>p(C_p(Y))=\aleph_0$ for
$Y=\beta\mathbb N$ (see [6], proof of Theorem 2.11).

\begin{corollary}
{\it Any Valdivia compact space $Y$ is metrizable if and only if
$p(C_p(Y))=\aleph_0$.}
\end{corollary}

\section{Linearly ordered compacta}

Let $(X,<)$ be a linearly ordered space. For every $x,y\in X$,
$x<y$, we put $(x,y)=\{z\in X: x<z<y\}$, $[x,y)=\{z\in X: x\leq
z<y\}$, $(x,y]=\{z\in X: x<z\leq y\}$ and $[x,y]=\{z\in X: x\leq
z\leq y\}$. Elements $x,y\in X$, $x<y$ such that $(x,y)=\O$ are
called {\it neighbor points}.

\begin{proposition}
{\it Let $(X,<)$ be a linearly ordered countably compact space and
let $A=\{(x,y)\in X^2: x<y \,\,\,and\,\,\, (x,y)=\O\}$. Then
$c(X)\leq p(C_p(X))$ and $|A|\leq p(C_p(X))$.}
\end{proposition}

{\bf Proof.} Let $B=\{(x_i,y_i)\in X^2:i\in I\}$ be a set which
satisfies the following conditions:

a) $x_i < y_i$ for every $i\in I$;

b) $y_i\leq x_j$ or $y_j\leq x_i$ for every distinct $i,j\in I$.

Put $U_i=\{f\in C_p(X): f(x_i)<0 \,\,\,and \,\,\, f(y_i)>1\}$ for
each $i\in I$. We show that $(U_i:i\in I)$ is a point-finite
family in $C_p(X)$. Suppose that $f\in C_p(X)$ such that a set
$J=\{i\in I: f\in U_i\}$ is infinite. Then in an infinite set
$\{x_i: i\in J\}$ there exists a strictly increasing or decreasing
sequence $(x_{i_n})^{\infty}_{n=1}$, $i_n\in J$. Since $X$ is a
countably compact space, there exists $z\in X$ such that $z=\lim
\limits_{n\to \infty}x_{i_n} = \lim \limits_{n\to \infty}y_{i_n}$.
Therefore $f(z)=\lim \limits_{n\to \infty}f(x_{i_n})\leq 0$ and
$f(z)=\lim \limits_{n\to \infty}f(y_{i_n})\geq 1$, which is
impossible.

Thus $(U_i:i\in I)$ is a point-finite family and $|I|\leq
p(C_p(X))$. Since any chain generates a set $B$ with $(x_i,y_i)\ne
\O$ for every $i\in I$, and $A$ generates a set $B$ with
$(x_i,y_i)= \O$ for every $i\in I$, we have the inequality
$c(X)\leq p(C_p(X))$ and $|A|\leq p(C_p(X))$.\hfill$\diamondsuit$

The following assertion seems to be known.

\begin{proposition}
{\it Let $(X,<)$ be a linearly ordered space. Then $d(X)\leq
c(X^2)$.}
\end{proposition}

{\bf Proof.} Let $A$ be a set of all isolated points of $X$ and
let $(W_i:i\in I)$ be a maximal family of nonempty sets
$W_i=(a_i,b_i)\times (c_i,d_i)$ in $X^2$ such that $b_i\leq c_i$
or $d_i\leq a_i$ for each $i\in I$. Clearly  $|A|\leq c(X^2)$ and
$|I|\leq c(X^2)$. Put $B=A\cup\{a_i, b_i, c_i, d_i: i\in I \}$.
Suppose that $U=X\setminus \overline {B}\ne \O$.

Choose a nonempty set $(x,y)\subseteq U$. Since $U\cap A=\O$ then
there exist $a,b,c,d\in (x,y)$ such that $a<b\leq c<d$ and sets
$(a,b)$ and $(c,d)$ are nonempty. Then an open in $X^2$ nonempty
set $W=(a,b)\times(c,d)$ such that $W\cap W_i =\O$ for every $i\in
I$. This contradicts the maximality of $(W_i:i\in
I)$.\hfill$\diamondsuit$

\begin{theorem}
{\it Let $(X, <)$ be an infinite linearly ordered compact. Then
$w(X)=p(C_p(X))$.}
\end{theorem}

{\bf Proof.} The inequality $p(C_p(X))\leq w(X)$ can be proved
like in the proof of Theorem 3.6.

Let $p(C_p(X))=\kappa$. Then a set $\tilde A=\{(x,y)\in X^2:
x<y\,\,\,and\,\,\,(x,y)=\O\}$ has the cardinality at most $\kappa$
by Proposition 4.1. Therefore a set $A=\{x,y\in X: (x,y)\in \tilde
A\}$ has the cardinality at most $\kappa$. Note that for arbitrary
dense in $X$ set $B$ the system $\{(x,y): x<y,\,\,x,y\in A\cup
B\}\bigcup \{[a,x): x\in A\cup B\}\bigcup\{(x,b]: x\in A\cup B\}$,
where $a=\min X$ and $b= \max X$, is a base of open in $X$ sets.
Therefore it is sufficient to prove that $d(X)\leq \kappa$.

Consider a system ${\mathcal V}$ of all open in $X$ nonempty set
$V$ with $d(V)\leq \kappa$. Choose a maximal disjoint system
${\mathcal U}\subseteq {\mathcal V}$. Proposition 4.1 implies
$|{\mathcal U}|\leq \kappa$. For every $U\in {\mathcal U}$ we
choose a set $B_U\subseteq U$ such that $|B_U|\leq \kappa$ and
$U\subseteq \overline{B_U}$. Put $B=A\bigcup
\bigcup\limits_{U\in{\mathcal U}}B_U$ and $Y=\overline{B}$. Then
$|B|\leq \kappa^2=\kappa$.

Suppose that $X\setminus Y\ne\O$. Choose $x_0, y_0\in X\setminus
Y$ such that $x_0<y_0$ and $X_0=[x_0,y_0]\subseteq X\setminus Y$.
We show that $c(X_0^2)\leq \kappa$.

Let $(W_i:i\in I)$ be a disjoint family of open sets in $X_0^2$.
Suppose that $|I|>\kappa$. For every $i\in I$ choose pair-wise
distinct points $a_i, b_i, c_i, d_i \in X_0$ such that $a_i<b_i$,
$c_i< d_i$ and $(a_i, b_i)\times (c_i, d_i)\subseteq W_i$ (it is
possible since $X_0$ has no neighbor point). Put $U_i=\{f\in
C_p(X_0): f(a_i)<0, f(b_i)>1, f(c_i)<0, f(d_i)>1\}$ for every
$i\in I$. Note that $p(C_p(X_0))\leq \kappa$, thus $(U_i:i\in I)$
is not point-finite in $C_p(X_0)$. There exist a function $f_0\in
C_p(X_0)$ and an infinite set $J\subseteq I$ such that $f_0\in
U_i$ for every $i\in J$. Since $f_0$ is continuous on linearly
ordered compact space $X_0$ which has no neighbor point, there
exists a finite set $K\subseteq X_0$ such that for every $x, y \in
X_0$ the inequality $|f_0(x) - f_0(y)|> 1$ implies $(x,y)\cap K
\ne \O$. Then $(a_i,b_i)\cap K\ne \O$, $(c_i,d_i)\cap K\ne \O$
therefore $W_i\cap K^2\ne \O$ for every $i\in J$ which contradicts
the disjointness of $(W_i: i\in I)$. Hence $|I|\leq \kappa$.

Thus $c(X_0^2)\leq \kappa$ and according to Proposition 4.2 we
obtain that $d(X_0)\leq \kappa$, which contradicts the maximality
of the system ${\mathcal U}$. Therefore $X\setminus Y =\O$,
$d(X)\leq \kappa$ and $w(X)\leq \kappa$.\hfill$\diamondsuit$

The following corollary gives a negative answer to Question 5.2
from  [7].

\begin{corollary}
{\it  A linearly ordered compact $X$ is metrizable if and only if
$p(C_p(X))=\aleph_0$.}
\end{corollary}


\begin{thebibliography}{00}

\bibitem [1] {1} {W.~Moran} {\it Separate continuity and support of
measures}, J. London. Math. Soc. {\bf 44} (1969), 320-324.

\bibitem [2] {2} {G.~Vera} {\it Baire measurability of separately continuous
functions}, Quart. J. Math. Oxford. {\bf 39, ¹ 153} (1988),
109-116.

\bibitem [3] {3} {V.V.~Mykhaylyuk, O.V.~Sobchuk} {\it Baire classification of
separately continuous functions and dependence upon countable
number of coordinates}, Nauk. Visn. Chern. Univ. {\bf Vyp.
191-192}, Matem., Chernivtsi: Ruta (2004), 116-118 (in Ukrainian).

\bibitem [4] {4} {V.V.~Tkachuk} {\it Cardinal invariants of the Suslin number
type}, Dokl. Akad. Nauk SSSR. {\bf 270, ¹ 4} (1983), 795-798 (in
Russian).

\bibitem [5] {5} {V.V.~Tkachuk} {\it A supertopological cardinal
invariants}, Vestnik Moskov. Univ. Ser. I Mat. Mekh. {\bf ¹ 4}
(1984), 26-29(in Russian).

\bibitem [6] {6} {A.V.~Arhangel'skii, V.V.~Tkachuk V.V.} {\it Calibers and
point-finite cellularity of the space $C_p(X)$ and some questions
of S.Gul'ko and M.Husek}, Topology Appl. {\bf 23, ¹ 1} (1986),
65-73.

\bibitem [7] {7} {O.G.~Okunev, V.V.~Tkachuk} {\it Lindel\"{o}f $\Sigma$-property
in $C_p(X)$ and $p(C_p(X))=\omega$ do not imply countable network
weight in $X$}, Acta Math. Hungar. {\bf 90, ¹ 1-2} (2001),
119-132.

\bibitem [8] {8} {N.D.~Kalamidas, G.D.~Spiliopoulos} {\it Compact sets in $C_p(X)$
and calibers}, Can. Math. Bull. {\bf 35, ¹ 4} (1992), 497-502.

\bibitem [9] {9} {V.V.~Mykhaylyuk} {\it One-point discontinuity set of
separately continuous functions on the product of two compact
spaces}, Ukr. Mat. Zhurn. {\bf 57, ¹ 1} (2005), 94-101 (in Ukrainian).


\bibitem [10] {10} {A.V.~Arhangelskii} {\it Topological spaces of functions},
Ì.: Izd-vo Moskovskogo Univ. (1989), 222 s.(in Russian).

\bibitem [11] {11} {R.~Engelking} {\it General topology},
M.: Mir (1986), 752 s. (in Russian).

\bibitem [12] {12} {M.~Talagrand} {\it Espaces de Banach faiblement ${\mathcal
K}$-analytiques}, Ann. of Math. {\bf 110} (1979), 407-438.


\end{thebibliography}
\end{document}